\documentclass{aims}
\usepackage{amsmath}
  \usepackage{paralist}
  \usepackage{graphics} 
  \usepackage{epsfig} 
\usepackage{graphicx}  \usepackage{epstopdf}
 \usepackage[colorlinks=true]{hyperref}
\hypersetup{urlcolor=blue, citecolor=red}

\usepackage[pagewise]{lineno}

\def\currentvolume{X}
 \def\currentissue{X}
  \def\currentyear{200X}
   \def\currentmonth{XX}
    \def\ppages{X--XX}
     \def\DOI{10.3934/xx.xx.xx.xx}

\renewcommand{\publname}{{\fontsize{7}{7}\selectfont Manuscript submitted to \hfill  {\href{http://dx.doi.org/\DOI}{doi:\DOI}}\\ Numerical Algebra, Control and Optimization}}

\def\currentvolume{} \def\currentissue{}  \def\currentyear{2018}   \def\currentmonth{}    \def\ppages{1--17}     \def\DOI{10.3934/xx.xx.xx.xx}

\newtheorem{theorem}{Theorem}[section]

\newtheorem{lemma}[theorem]{Lemma}

\theoremstyle{definition}
\newtheorem{definition}[theorem]{Definition}
\newtheorem{remark}{Remark}

\newcommand{\R}{\mathbb{R}}
\newcommand{\N}{\mathbb{N}}
\newtheorem{example}{Example}

\title[Upper and Lower Bounds for the HOMO-LUMO Spectral Gap] 
      {On Construction of Upper and Lower Bounds for the HOMO-LUMO Spectral Gap}

\author[So\v{n}a Pavl\'{\i}kov\'a, Daniel \v{S}ev\v{c}ovi\v{c}]{}


 \email{sona.pavlikova@stuba.sk}
 \email{sevcovic@fmph.uniba.sk}

\thanks{The authors are supported by VEGA grant 1/0062/18}

\thanks{$^*$ Corresponding author:  Daniel \v{S}ev\v{c}ovi\v{c}}

\begin{document}

\maketitle
\centerline{\scshape So\v{n}a Pavl\'{\i}kov\'a}
\medskip
{\footnotesize
 \centerline{
 Institute of Information Engineering, Automation, and Mathematics}
   \centerline{FCFT, Slovak Technical University}
   \centerline{812 37 Bratislava, Slovakia}
} 

\medskip

\centerline{\scshape Daniel \v{S}ev\v{c}ovi\v{c}$^*$}
\medskip
{\footnotesize
 \centerline{Department of Applied Mathematics and Statistics}
 \centerline{FMFI, Comenius University}
 \centerline{842 48 Bratislava, Slovakia}
}

\bigskip

 \centerline{(Communicated by the associate editor name)}

\begin{abstract}
In this paper we study spectral properties of graphs which are constructed from two given invertible graphs by bridging them over a bipartite graph. We analyze 
the so-called HOMO-LUMO spectral gap which is the difference between the smallest positive and largest negative eigenvalue of the adjacency matrix of a graph. We investigate  its dependence on the bridging bipartite graph and we construct a mixed integer semidefinite program for maximization of the HOMO-LUMO gap with respect to the bridging bipartite graph. We also derive upper and lower bounds for the optimal HOMO-LUMO spectral graph by means of semidefinite relaxation techniques. Several computational examples are also presented in this paper. 
\end{abstract}

\section{Introduction} \label{sec-intro}

The spectrum $\sigma(G)$ of an undirected graph $G$ consists of eigenvalues of its adjacency symmetric $n\times n$ matrix ${\mathcal A}(G)$, i.e. $\sigma(G)=\{\lambda_k(G), k=1,\cdots, n, \ \lambda_k(G)$ is an eigenvalue of  ${\mathcal A}(G)\}$, where 
$\lambda_1(G) \ge \cdots \ge \lambda_n(G)$ (cf. \cite{Cvetkovic1988, Brouwer2012}). If the spectrum does not contain zero there exists the inverse matrix $A^{-1}$ of the adjacency matrix $A={\mathcal A}(G)$, and the graph $G_A$ is called invertible. 

The concept of an inverse graph has been introduced by Godsil \cite{Godsil1985}. In addition to invertibility of the adjacency matrix it is required that $A^{-1}$ is diagonally similar to a nonnegative or nonpositive integral matrix (cf. Godsil \cite{Godsil1985}, Pavl\'{\i}kov\'a and \v{S}ev\v{c}ovi\v{c} \cite{Pavlikova2016}). Notice that the least positive eigenvalue of a graph is the reciprocal value of the maximal eigenvalue of the inverse graph. Therefore properties of inverse graphs can be used in estimation of the least positive eigenvalue (cf. Pavl\'{\i}kov\'a et al. \cite{Pavlikova1990, Pavlikova2015, Pavlikova2016}). 

In many applied fields, e.g. theoretical chemistry, biology, or statistics, spectral indices and properties of graphs representing structure of chemical molecules or transition diagrams for finite Markov chains play an important role (cf. Cvetkovi\'c \cite{Cvetkovic1988,Cvetkovic2004}, Brouwer and Haemers \cite{Brouwer2012} and references therein). In the last decades, various graph energies and indices have been proposed and analyzed. For instance, the sum of absolute values of eigenvalues is referred to as the matching energy index (cf. Chen and Liu \cite{Lin2016}), the maximum of the absolute values of the least positive and largest negative eigenvalue is known as the HOMO-LUMO index (see Mohar \cite{Mohar2013,Mohar2015}, Li \emph{et al.} \cite{Li2013}, Jakli\'c  \emph{et al.} \cite{Jaklic2012}, Fowler \emph{et al.} \cite{Fowler2010}), their difference is the  HOMO-LUMO separation gap (cf.  Gutman and Rouvray \cite{Gutman1979}, Li \emph{et al.} \cite{Li2013}, Zhang and An \cite{Zhang2002}, Fowler  \emph{et al.} \cite{Fowler2001}).

In computational chemistry, eigenvalues of a graph describing an organic molecule are related to energies of molecular orbitals. Following H\"uckel's molecular orbital method \cite{Huckel1931}  (see also Pavl\'{\i}kov\'a and \v{S}ev\v{c}ovi\v{c} \cite{Pavlikova2016-CMMS}), the energies $E_k, k=1,\cdots, n$, are the eigenvalues of the Hamiltonian matrix $H$ and its eigenvectors are orbitals. The square symmetric matrix $H$ has the following elements:
\begin{itemize}
 \item[] $H_{ii} = \alpha$ for the carbon C atom at the $i$-th vertex, and $H_{ii} = \alpha + h_A\beta$ for other atoms A, where $\alpha<0$ is the Coulomb integral and $\beta<0$ is the resonance integral;
\item[] $H_{ij} = \beta$ if both vertices $i$ and $j$ are carbon C atoms,  $H_{ij} = k_{AB}\beta$ for other neighboring atoms A and B;
\item[] $H_{ij} = 0$ otherwise.
\end{itemize}
The atomic constants $h_A, k_{AB}$ have to be specified ($h_C=k_{CC}=0$). For instance, the molecule of pyridine contains one atom of nitrate N and five atoms of carbon C. Clearly, in the case of pure hydrocarbon we have $H = \alpha I + \beta A$ where $I$ is the identity and $A$ is the adjacency matrix of the molecular structural graph $G$. Hence $E_k=\alpha +\beta \lambda_k$. Now, the energy $E_{HOMO}$ of the highest occupied molecular orbital (HOMO) corresponds to the eigenvalue $\lambda_{HOMO}=\lambda_k$ where $k=n/2$ for $n$ even and $k=(n+1)/2$ for $n$ odd. The energy $E_{LUMO}$ of the lowest unoccupied molecular orbital (LUMO) corresponds to the subsequent eigenvalue $\lambda_{LUMO}=\lambda_{k+1}$ for $n$ even, and $\lambda_{LUMO}=\lambda_{k}$ for $n$ odd. The HOMO-LUMO separation gap is the difference between $E_{LUMO}$ and $E_{HOMO}$ energies, i.e.  $E_{LUMO} - E_{HOMO} = -\beta ( \lambda_{HOMO} - \lambda_{LUMO}) \ge 0$ because $\beta<0$. The so-called properly closed shells have the property $\lambda_{HUMO}>0>\lambda_{LUMO}$ containing either zero or two electrons are called closed shells for which $n$ is even (cf. Fowler and Pisanski \cite{Fowler2010}). For such orbital systems, the HOMO-LUMO separation gap is equal to the energy difference $E_{LUMO} - E_{HOMO} = -\beta \Lambda_{HL}(G_A)$
where 
\begin{equation}
\Lambda_{HL}(G_A) = \check{\lambda}^+(G_A) - \hat{\lambda}^-(G_A).
\label{HLgap}
\end{equation}
Here $\check{\lambda}^+(G_A)=\lambda_k$ is the smallest positive eigenvalue, and  $\hat{\lambda}^-(G_A)=\lambda_{k+1}$  is the largest negative eigenvalue  of the adjacency matrix $A$ of the structural molecular graph $G_A$ (cf. \cite{Fowler2010}). According to Aihara \cite{Aihara1999JCP,Aihara1999TCH} the large HOMO-LUMO gap implies high kinetic stability and low chemical reactivity of the molecule, because it is energetically unfavorable to add electrons to a high-lying LUMO orbital. Notice that the HOMO-LUMO energy gap is generally decreasing with the size $n$ of the structural graph (cf. Bacalis and Zdetsis \cite{Bacalis2009}).

In this paper, our goal is to investigate extremal properties of the HOMO-LUMO spectral gap $\Lambda_{HL}(G_A)$. We show how to represent $\Lambda_{HL}(G_A)$ by means of the optimal solution to a convex semidefinite programming problem (Section 2). We study spectral properties of graphs which can be constructed from two given (not necessarily bipartite) graphs by bridging them over a bipartite graph (Section 3). We analyze their HOMO-LUMO spectral gap of such a bridged graph and its dependence on the bridging bipartite graph. Finding an optimal bridging bipartite graph leads to a mixed integer nonconvex optimization problem with linear matrix inequality constraints (Section 4). We prove that the optimal HOMO-LUMO spectral gap can be obtained by solving a mixed integer semidefinite convex program. The optimization problem is, in general, NP hard  (Section 5). This is why we also derive upper (Section 6) and lower (Section 7) bounds for the optimal HOMO-LUMO spectral graph by means of semidefinite relaxation techniques which can be solved in a fast and computationally efficient way. Various computational examples of construction of the optimal bridging graph are presented in Section 8.

\section{Semidefinite programming representation of the HOMO-LUMO spectral gap}

The HOMO-LUMO spectral gap of a graph $G_C$ is defined as follows:
\[
 \Lambda_{HL}(G_C) = \check{\lambda}^+(G_C) - \hat{\lambda}^-(G_C),
\]
where $\check{\lambda}^+(G_C) \ge 0$ is the smallest nonnegative eigenvalue, and $\hat{\lambda}^-(G_C)\le 0$ is the largest nonpositive eigenvalue  of the adjacency matrix $C$. Notice that the spectrum $\sigma(G_C) = \sigma(C)$ of a nontrivial graph $G_C$ without loops must contain negative as well as positive eigenvalues because the trace $Tr(C) =\sum_{\lambda\in\sigma(C)} \lambda = 0$. 
Clearly, if the graph $G_C$ is invertible then $\check{\lambda}^+(G_C) > 0$ and $\hat{\lambda}^-(G_C)< 0$ and so $\Lambda_{HL}(G_C)>0$, otherwise $\Lambda_{HL}(G_C)=0$.

\subsection{Semidefinite representation of the HOMO-LUMO gap}

Suppose that a graph $G_C$ is invertible. Following \cite{Pavlikova2016} the smallest positive and largest negative eigenvalues of $G_C$ can be expressed as follows:
\[
  \check{\lambda}^+(G_C) = \frac{1}{\lambda_{max}(C^{-1})}, \qquad \hat{\lambda}^-(G_C) = \frac{1}{\lambda_{min}(C^{-1})},
\]
where $\lambda_{max}(C^{-1})>0$ and $\lambda_{min}(C^{-1})= - \lambda_{max}(-C^{-1})<0$ are the maximum and minimum eigenvalues of the inverse matrix $C^{-1}$, respectively. We denote by $\preceq$ the L\"owner partial ordering on symmetric matrices, i.e. $A\preceq B$ iff the matrix $B-A$ is a positive semidefinite matrix, that is $B-A\succeq 0$. The maximal and minimal eigenvalues of $C^{-1}$ can be expressed as follows:
\[
0< \lambda_{max}(C^{-1}) = \min_{C^{-1}\preceq t I} t, \qquad 0 > \lambda_{min}(C^{-1}) = \max_{ s I \preceq C^{-1}} s,
\]
(see e.g. \cite{bova}, \cite{Cvetkovic2004}). 
Since $\{ t,\ C^{-1} \preceq t I \} \subset (0,\infty)$ and $\{ s,\ s I \preceq C^{-1} \} \subset (-\infty, 0)$ then, by using the substitution $\mu=1/t, \eta=-1/s$, we obtain the following characterization of the lowest positive and largest negative eigenvalues of the graph $G_C$:
\begin{equation}
 \check{\lambda}^+(G_C) = \max_{\mu C^{-1}\preceq I} \mu, \qquad \hat{\lambda}^-(G_C) = - \max_{- \eta C^{-1}\preceq I} \eta.
\label{lambdapm}
\end{equation}

As a consequence, we obtain the following semidefinite representation of the HOMO-LUMO spectral gap for a vertex labeled invertible graph $G_C$ without loops. Then the HOMO-LUMO spectral gap $\Lambda_{HL}(G_C)$ of the graph $G_C$ is the optimal value of the following semidefinite programming problem:
\begin{eqnarray}
\Lambda_{HL}(G_C) &=& \max_{\mu,\eta\ge0}  \quad \mu+\eta
\\
&& s.t. \quad  \mu C^{-1} \preceq I,  \nonumber
\\
&& \quad\ \   -\eta  C^{-1} \preceq I. \nonumber
\label{homolumo}
\end{eqnarray}
(cf. Pavl\'{\i}kov\'a and \v{S}ev\v{c}ovi\v{c} \cite{Pavlikova2016-CMMS}).

\section{Graphs bridged over a bipartite graph}

In this section we introduce a notion of a graph which is constructed from two given graphs $G_A$ and $G_B$ by bridging vertices of $G_A$ to vertices of $G_B$.  More, precisely, let $G_A$ and $G_B$ be two undirected vertex-labeled graphs on $n$ and $m$ vertices without loops, respectively. In general, we do not assume that $G_A$ and $G_B$ are bipartite graphs. Let $G_K$ be a $(n,m)$-bipartite graph on $n+m$ vertices with the  adjacency matrix:
\begin{equation}
{\mathcal A}(G_K) = \left( 
\begin{array}{cc}
0 & K\\
K^T & 0
\end{array}
\right),
\label{bipartite}
\end{equation}
where $K$ is an $n\times m$ matrix containing $\{0,1\}$-elements only.

\begin{figure}[htp]
\begin{center}
\includegraphics[width=7truecm]{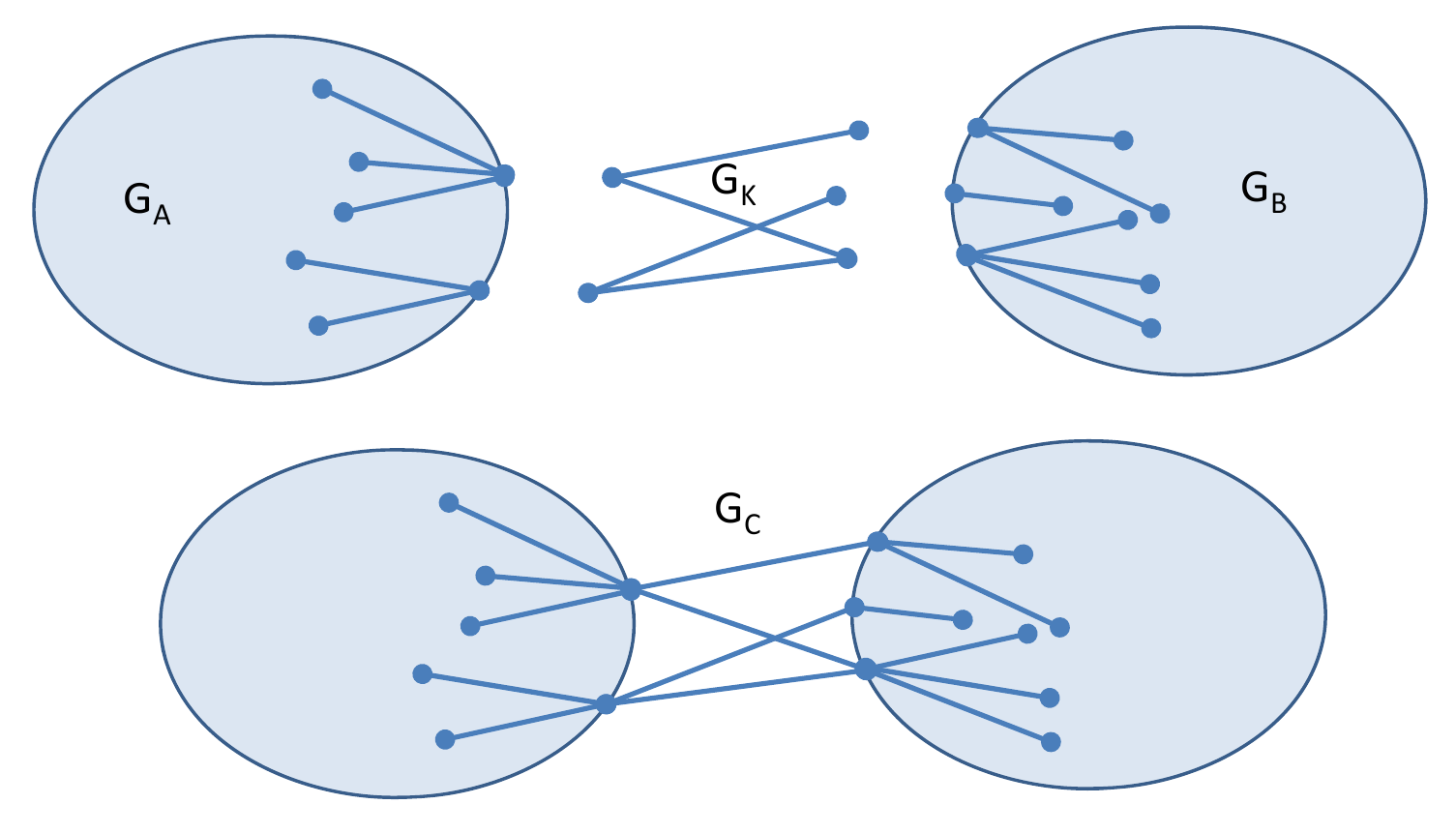}
\end{center}
\caption{
A bridged graph $G_C={\mathcal B}_K(G_A,G_B)$ through a bipartite graph $G_K$.}
\label{fig-bipartite}
\end{figure}

By ${\mathcal B}_K(G_A,G_B)$ we shall denote the graph $G_C$ on $n+m$ vertices which is obtained by bridging the vertices of the graph $G_A$ to the vertices of $G_B$ through the $(n,m)$-bipartite graph $G_K$, i.e. its adjacency matrix $C={\mathcal A}(G_C)$ of the graph $G_C$ has the form:
\begin{equation}
C  = \left( 
\begin{array}{cc}
A & K\\
K^T & B
\end{array}
\right),
\label{matrixC}
\end{equation}

In what follows, we will assume that the adjacency matrices $A$ and $B$ are symmetric $n\times n$ and $m\times m$ invertible matrices, respectively.

\begin{theorem}\label{theo-1}
Let $G_A$ and $G_B$ be two undirected vertex-labeled invertible graphs on $n$ and $m$ vertices, respectively. Let $G_K$ be a $(n,m)$-bipartite graph.  Let $G_C={\mathcal B}_K(G_A,G_B)$ be the graph which is constructed by bridging the graphs $G_A$ and $G_B$ through the bipartite graph $G_K$. 

Then the graph $G_C$ is invertible if and only if the $n\times n$ matrix  $S= A- K B^{-1} K^T$ is invertible. In this case we have 
\begin{eqnarray}
C^{-1}&=&\left( 
\begin{array}{cc}
A & K\\
K^T & B
\end{array}
\right)^{-1}
=
\left( 
\begin{array}{cc}
S^{-1} & - S^{-1} K B^{-1} \\
- B^{-1} K^T S^{-1} & B^{-1} + B^{-1} K^T S^{-1} K B^{-1}
\end{array}
\right).
\nonumber 
\\
&=&
Q^T \left( 
\begin{array}{cc}
S^{-1} & 0 \\
0 & B^{-1}
\end{array}
\right)
Q,
\label{invC2}
\end{eqnarray}
where $Q$ is an invertible matrix with the inverse $Z=Q^{-1}$ given by:
\[
Q =
\left( 
\begin{array}{cc}
I & - K B^{-1}  \\
0 &  I
\end{array}
\right), \qquad Z =
\left( 
\begin{array}{cc}
I & K B^{-1} \\
0 &  I
\end{array}
\right).
\]

\end{theorem}

\noindent P r o o f.  The proof is a direct consequence of the Schur complement theorem (see e.~g. \cite[Theorem A.6]{Maja2013}). Indeed, $C\left(\begin{array}{c} x\\ y\end{array}\right) = \left(\begin{array}{c} 0\\ 0\end{array}\right) $ if and only if $Ax + Ky =0$ and $K^Tx + By=0$, that is, $S x = (A- K B^{-1} K^T)x =0$. As $x\not=0\Leftrightarrow y\not=0$ we have $C$ is invertible if and only if $S$ is invertible. The rest of the proof is a straightforward verification of the form of the inverse matrix $C^{-1}$. 
\hfill $\diamondsuit$

\subsection{Semidefinite representation of the HOMO-LUMO gap for a bridged graph}

Now, let $G_C = {\mathcal B}_K(G_A,G_B)$ be the graph obtained from graphs $G_A$ and $G_B$ by bridging them through a bipartite graph $G_K$ with adjacency matrix  $K$ (\ref{bipartite}). 

Then, for any $\mu\ge 0$, we have $\mu C^{-1} \preceq I$ if and only if $\mu Z^T C^{-1} Z \preceq Z^T Z$, i.e., 
\[
 \mu \left( 
\begin{array}{cc}
S^{-1} & 0 \\
0 &  B^{-1}
\end{array}
\right) 
\preceq Z^T Z = 
\left( 
\begin{array}{cc}
I  &  K B^{-1} \\
B^{-1} K^T &  I + B^{-1} K^T K B^{-1}
\end{array}
\right).
\]
Therefore,
\begin{equation}
\mu C^{-1} \preceq I \quad\Leftrightarrow \quad 
\left( 
\begin{array}{cc}
I - \mu S^{-1} &  K B^{-1} \\
B^{-1} K^T &  I - \mu B^{-1} + B^{-1} K^T K B^{-1}
\end{array}
\right) \succeq 0.
\label{ineqmu}
\end{equation}
Similarly, 
\begin{equation}
-\eta C^{-1} \preceq I \quad\Leftrightarrow \quad 
\left( 
\begin{array}{cc}
I + \eta S^{-1} &  K B^{-1} \\
B^{-1} K^T &  I + \eta B^{-1} + B^{-1} K^T K B^{-1}
\end{array}
\right) \succeq 0.
\label{ineqeta}
\end{equation}
With regard to (\ref{homolumo}) we obtain the following representation of the HOMO-LUMO spectral gap $\Lambda_{HL}(G_C)$ for a the bridged graph:
\begin{eqnarray}
\label{homolumobridged-S}
\Lambda_{HL}(G_C) &=& \max_{\scriptsize
\begin{array}{c}
\mu, \eta \ge 0
\end{array} 
}  
\quad \mu+\eta
\\
s.t.&& 
\left( 
\begin{array}{cc}
I - \mu S^{-1}  &  K   B^{-1} \\
B^{-1} K^T   &  I - \mu B^{-1} + B^{-1} K^T K  B^{-1}
\end{array}
\right) \succeq 0,
\nonumber
\\
&& 
\left( 
\begin{array}{cc}
I +\eta S^{-1} &  K  B^{-1} \\
B^{-1} K^T  &  I + \eta B^{-1} + B^{-1} K^T K  B^{-1}
\end{array}
\right) \succeq 0. 
\nonumber
\end{eqnarray}

Since for the Schur complement we have $S= A- K B^{-1} K^T$ then the matrix inequality constraints appearing in (\ref{homolumobridged-S}) represent, in general, nonconvex constraints with respect to the matrix $K$. To overcome this difficulty we further restrict the class of bipartite graphs $G_K$ bridging $G_A$ to $G_B$ to those turning  (\ref{homolumobridged-S}) to a convex semidefinite program in the $K$ variable.

\begin{definition}\cite{Pavlikova2016} 
\label{def-arbitrarily}
Let $G_B$ be an undirected vertex-labeled graph on $m$ vertices with an invertible adjacency matrix $B$. We say that $G_B$ is arbitrarily bridgeable over the first $\{1, \cdots, k_B\}$ vertices of $G_B$ if the $k_B\times k_B$ upper principal sub-matrix of $B^{-1}$ is a null matrix, i.e.  $E B^{-1} E^T=0$ where $E=(I, 0)$ is a $k_B\times m$ block matrix and $I$ is a $k_B\times k_B$ identity matrix. 

A graph $G_B$ is said to be arbitrarily bridgeable over the subset $\{ i_1, \cdots, i_{k_B}\}$ of vertices of $G_B$ if there exists a permutation $P$ of its vertices such that $i_k\mapsto k, k=1,\cdots, k_B,$ and $E \tilde B^{-1} E^T=0$ where $\tilde B= P^T B P$. 
\end{definition}

Notice that if $G_B$ is arbitrarily bridgeable then $k_B\le m/2$ because there is no regular $m\times m$  matrix $B^{-1}$ such that $E B^{-1} E^T =0$ for $k_B > m/2$.

Using the notion of an arbirtarily bridgeable graph we conclude the following theorem:

\begin{theorem}\label{th-semidefinitebridged}
Let $G_A$ and $G_B$ be undirected vertex-labeled invertible graphs on $n$ and $m$ vertices without loops, respectively. Assume that $G_B$ is arbitrarily bridgeable over the first $\{1, \cdots, k_B\}$ vertices of $G_B$. If the $n\times m$ matrix $K$ has zero last $m-k_B$ columns, i.e. $K_{ij}=0$ for $j=k_B+1, \cdots m$, then $K B^{-1} K^T = 0$, and, consequently, for the Schur complement $S$ we have   $S= A- K B^{-1} K^T = A$, and $S^{-1}= A^{-1}$. 

Moreover, the HOMO-LUMO spectral gap $\Lambda_{HL}(G_C)$ for the bridged graph $G_C={\mathcal B}_K(G_A,G_B)$ through the bipartite graph $G_K$ is the optimal value of the following semidefinite programming problem:
\begin{eqnarray}\small
\label{homolumobridged}
\Lambda_{HL}(G_C) &=& \max_{\scriptsize
\begin{array}{c}
\mu, \eta \ge 0
\end{array} 
}  
\quad \mu+\eta
\\
s.t.&& 
\left( 
\begin{array}{cc}
I - \mu A^{-1}  &  K  B^{-1} \\
B^{-1} K^T   &  I - \mu B^{-1} + B^{-1}  K^T K B^{-1}
\end{array}
\right) \succeq 0,
\nonumber
\\
&& 
\left( 
\begin{array}{cc}
I +\eta A^{-1} &  K  B^{-1} \\
B^{-1} K^T  &  I + \eta B^{-1} + B^{-1}  K^T K  B^{-1}
\end{array}
\right) \succeq 0.
\nonumber
\end{eqnarray}
\end{theorem}

\section{Construction of an optimal bridging bipartite graph by means of a mixed integer nonlinear programming problem}

In this section we focus our attention on extremal properties of the HOMO-LUMO spectral gap for bridged graphs. Given an invertible graph $G_A$ and arbitrarily bridgeable invertible graph $G_B$, over the first $\{1, \cdots, k_B\}$ vertices of $G_B$, our goal is to find an optimal bridging graph $G_K$ (see (\ref{bipartite})) such that $K_{ij}=0$ for $j=k_B+1, \cdots, m$ and the HOMO-LUMO spectral gap $\Lambda_{HL}(G_C)$ is maximal, where $G_C={\mathcal B}_K(G_A, G_B)$.

Using representation of $\Lambda_{HL}(G_C)$ for the graph $G_C={\mathcal B}_K(G_A,G_B)$ (see Theorem~\ref{th-semidefinitebridged}), the maximal HOMO-LUMO gap $\Lambda^{opt}_{HL}=\Lambda^{opt}_{HL}(G_A, G_B)$ with respect to a bipartite matrix $K$ is given as the optimal value of the following mixed integer nonlinear optimization problem:

\begin{eqnarray}\small
\label{homolumoopt}
\Lambda^{opt}_{HL} &=& \max_{\scriptsize
\begin{array}{c}
\mu, \eta \ge 0 \\
K, W 
\end{array} 
}  
\quad \mu+\eta
\\
s.t.&& 
\left( 
\begin{array}{cc}
I - \mu A^{-1}  &  K  B^{-1} \\
B^{-1}  K^T   &  I - \mu B^{-1} + B^{-1} W  B^{-1}
\end{array}
\right) \succeq 0,
\nonumber
\\
&& 
\left( 
\begin{array}{cc}
I +\eta A^{-1} &  K   B^{-1} \\
B^{-1} K^T  &  I + \eta B^{-1} + B^{-1} W  B^{-1}
\end{array}
\right) \succeq 0,
\nonumber
\\
&& \ \ 
W = K^T K, \quad K_{ij}\in \{0,1\}\quad\hbox{for all}\  i,j, \quad
\sum_{k,l} K_{kl} \ge 1,
\nonumber \\
&& K_{ij}=0 \ \ \hbox{for} \ j=k_B+1, \cdots, m,\ \  i=1,\cdots, n.
\nonumber
\end{eqnarray}
Notice that the condition $K\not=0$ for a binary matrix $K$ is equivalent to the condition $\sum_{k,l} K_{kl} \ge 1$. The  objective function as well as the first two matrix inequality constraints in the optimization problem (\ref{homolumoopt}) are linear\footnote{Convex semidefinite problems with linear matrix inequality constraints can be solved by means of computational Matlab toolboxes available for semidefinite programming, e.g. SeDuMi solver developed by J.~Sturm \cite{sturm} with Yalmip Matlab programming framework due to J. L\"ofberg \cite{Lofberg2004}.} in the variables $\mu,\eta, K, W$. However, the last two constraints in (\ref{homolumoopt})  make the problem considerably harder to solve because of the nonconvex constraint $W = K^T K$ and the binary constraint $K_{ij}\in \{0,1\}$. It  means that (\ref{homolumoopt}) is a mixed integer nonconvex programming problem which is, in general, NP-hard to solve.

\section{Construction of upper bounds for the HOMO-LUMO spectral gap by semidefinite relaxation techniques}

In the field of solving mixed integer nonconvex problems  various techniques have been developed in the last decades. We refer the reader to the book \cite{bova} by Boyd and Vanderberghe on recent developments on semidefinite relaxation methods for solving nonconvex and mixed integer nonlinear optimization problems. In general, semidefinite relaxations of an original nonconvex problem can be constructed by means of the second Lagrangian dual problem which is already a convex semidefinite problem (see e.g. \v{S}ev\v{c}ovi\v{c} and Trnovsk\'a \cite{ST}). 

\subsection{Mixed semidefinite-integer relaxation}
In order to construct a suitable convex programming relaxation of (\ref{homolumoopt}) we have to enlarge the domain of variables $\mu,\eta, K, W$.   Notice that the integer constraint $K_{ij}\in \{0,1\}$ is equivalent to the equality: $K_{ij}=K^2_{ij}$. Moreover, from the constraint $W=K^T K$ we deduce $W_{ij}\in \N^+_0$ and $W_{jj} = \sum_{l} K^2_{lj} = \sum_{l} K_{lj}$. The nonconvex constraint $W=K^T K$ can be relaxed by a convex matrix inequality constraint $W\succeq K^T K$. Using the Schur complement theorem (cf. \cite{Maja2013}), it can be rewritten as a linear matrix inequality constraint:
\[
 W\succeq K^T K \quad\Leftrightarrow\quad 
\left( \begin{array}{cc}
W  &  K^T \\
K  &  I
\end{array}\right)\succeq 0.
\]
Hence the nonconvex-integer programming problem (\ref{homolumoopt}) can be relaxed by means of the following mixed integer semidefinite programming problem with linear matrix inequality constraints and integer constraints for the upper bound approximation $\overline{\Lambda}^{sir}_{HL}=\overline{\Lambda}^{sir}_{HL}(G_A,G_B)$:

\begin{eqnarray}
\label{homolumosir}
\overline{\Lambda}^{sir}_{HL} &=& \max_{\scriptsize
\begin{array}{c}
\mu, \eta \ge 0 \\
K, W 
\end{array} 
}  
\quad \mu+\eta
\nonumber 
\\
s.t.&& 
\left( 
\begin{array}{cc}
I - \mu A^{-1}  &  K   B^{-1} \\
B^{-1}  K^T  &  I - \mu B^{-1} + B^{-1} W  B^{-1}
\end{array}
\right) \succeq 0,
\nonumber
\\
&& 
\left( 
\begin{array}{cc}
I +\eta A^{-1} &  K  B^{-1} \\
B^{-1} K^T   &  I + \eta B^{-1} + B^{-1} W B^{-1}
\end{array}
\right) \succeq 0,
\\
&&
\left( \begin{array}{cc}
W  &  K^T \\
K  &  I
\end{array}\right)\succeq 0,\nonumber
\\
&&
K_{ij}\in \{0,1\}, 
\ \  W_{ij}\in\N^+_0, \ \ W_{jj} = \sum_{l} K_{lj}
\quad\hbox{for all}\ i,j, \ \sum_{k,l} K_{kl} \ge 1.  \nonumber
\\
&& K_{ij}=0 \ \ \hbox{for} \ j=k_B+1, \cdots, m, \ \ i=1,\cdots, n.
\nonumber
\end{eqnarray}

It is worth noting that if $(\hat\mu, \hat\eta, \hat K, \hat W)$ is the optimal solution to the mixed integer semidefinite programming problem (\ref{homolumosir}) then $(\hat\mu, \hat\eta, \hat K, \hat W)$ is also feasible for (\ref{homolumoopt}) because $\hat W=\hat K^T \hat K$. Indeed, if we denote $L=\hat W - \hat K^T \hat K$ then $L\succeq 0$ and $L_{jj}=\hat W_{jj} - \sum_l \hat K_{lj}^2 = \hat W_{jj} - \sum_l \hat K_{lj} =0$. Hence $diag(L)=0$ and so $L=0$, as claimed. Consequently, the HOMO-LUMO gap $\Lambda_{HL}({\mathcal B}_{\hat K}(G_A,G_B) ) = \overline{\Lambda}^{sir}_{HL}(G_A,G_B)$. Hence
\[
 \Lambda^{opt}_{HL}(G_A,G_B) = \overline{\Lambda}^{sir}_{HL}(G_A,G_B).
\]

Next we present a sample code for solving the mixed integer semidefinite programming problem (\ref{homolumosir}) for construction of the optimal bridging for maximal HOMO-LUMO spectral gap $\overline{\Lambda}^{sir}_{HL}(G_A, G_B) = \overline{\Lambda}^{opt}_{HL}(G_A, G_B)$. We employed the Matlab programming environment Yalmip which is capable of solving mixed integer problems with semidefinite linear matrix inequality constraints due to L\"ofberg \cite{Lofberg2004}). The structure of the code is shown in Table~\ref{tab-code}. After declaring classes of variables and setting the constraints, then the main solver routine {\tt solvsdp} is executed. It is designed for solving minimization problem. It employs SeDuMi semidefinite programming solver (cf. Sturm \cite{sturm}) as the lower solver and branch and bound integer rounding solver as the upper solver. 

\begin{table}[htp]
\label{tab-code}
\caption{A sample Matlab code for computing mixed integer semidefinite programming problem (\ref{homolumosir}). The output of the program is the optimal value 
$\Lambda^{opt}_{HL}(G_A,G_B) = \overline{\Lambda}^{sir}_{HL}(G_A,G_B)$.
}

\begin{center}
\hrule
\vskip 0.5truemm
\hrule

{\footnotesize
\begin{verbatim}

mu=sdpvar(1); eta=sdpvar(1); W=intvar(m,m); K=binvar(n,m);
ops=sdpsettings('solver','bnb','bnb.maxiter', bnbmaxiter);

Fconstraints=[...
    [[W, K'];
    [K, eye(n,n)]
    ]>=0, ...
    mu>=0, eta>=0, ...
    [[eye(n,n) - mu*inv(A), K*inv(B)];
    [inv(B)*K', eye(m,m) - mu*inv(B) + inv(B)*W*inv(B)]
    ] >= 0, ...
    [[eye(n,n) + eta*inv(A), K*inv(B)];
        [inv(B)*K', eye(m,m) + eta*inv(B) + inv(B)*W*inv(B)]
    ] >= 0, ...
    sum(K(:,:))==diag(W)', sum(K(:))>=1, ...
    vec(W(:))>=0, 0<=vec(K(:))<=1, ...
    sum([[A, K]; [K', B] ])<=maxdegree*ones(1,n+m), ...,
    K*[zeros(kB,m-kB); eye(m-kB,m-kB)] == zeros(n, m-kB), ...
    ];

solvesdp(Fconstraints, -mu-eta, ops)

LambdaSIR = double(mu + eta)

\end{verbatim}
}

\hrule

\end{center}

\end{table}

\subsection{Full semidefinite relaxation}

Next, we further relax the binary and integer constraints appearing in (\ref{homolumosir}). The integer constraint $K_{ij}\in\{0,1\}$ can be relaxed by the box convex inequality constraints: $0\le K_{ij}\le 1$ for all $i,j$. Clearly, such a relaxation may lead to a non-integer optimal matrix $K$. The maximization problem for the full semidefinite relaxation of the HOMO-LUMO spectral gap $\overline{\Lambda}_{HL}^{sdp}=\overline{\Lambda}_{HL}^{sdp}(G_A,G_B)$ can be formulated as follows:
\begin{eqnarray}
\label{homolumosdp}
\overline{\Lambda}_{HL}^{sdp} &=& \max_{\scriptsize
\begin{array}{c}
\mu, \eta \ge 0 \\
K, W
\end{array} 
}  
\quad \mu+\eta
\nonumber \\
s.t.&& 
\left( 
\begin{array}{cc}
I - \mu A^{-1}  & K  B^{-1} \\
B^{-1} K^T   &  I - \mu B^{-1} + B^{-1} W  B^{-1}
\end{array}
\right) \succeq 0,
\nonumber
\\
&& 
\left( 
\begin{array}{cc}
I +\eta A^{-1} &  K   B^{-1} \\
B^{-1}  K^T  &  I + \eta B^{-1} + B^{-1} W B^{-1}
\end{array}
\right) \succeq 0,
\\
&&
\left( 
\begin{array}{cc}
W    & K^T        \\
K    &  I \\ 
\end{array} 
\right)\succeq 0, \nonumber 
\\
&& 
0\le K_{ij} \le 1,\ \  W_{jj} = \sum_{l} K_{lj}, \ \ 
W_{ij}\ge 0\  \hbox{for all}\ i,j,\ \ \sum_{k,l} K_{kl}\ge 1, 
\nonumber
\\
&& K_{ij}=0 \ \ \hbox{for} \ j=k_B+1, \cdots, m, \ \ i=1,\cdots, n.
\end{eqnarray}

In order to compute $\overline{\Lambda}^{sdp}_{HL}(G_A, G_B)$ the full semidefinite relaxation (\ref{homolumosdp})  we have to change the specification of real variables, i.e. {\tt W=sdpvar(m,m); K=sdpvar(n,m)} and add the box constraint {\tt  0<=vec(K(:))<=1} in the  code shown in Table~\ref{tab-code}. 

\begin{remark}
Following the recent paper by Kim, Kojima and Toh \cite{KKT2016} the  box constraint $0\le K_{ij}\le 1$ can be further enhanced by introducing a slack variable $\tilde K$ where $\tilde K_{ij}= 1-K_{ij}$. Then $K_{ij}\in\{0,1\}$ if and only if $K_{ij}\tilde K_{ij}=0$ for all $i,j$. It is equivalent to the condition $V_{jj}=0$ for each $j$, where $V = \tilde K^T K$. Next, the nonconvex matrix constraints $W= K^T K, \tilde W = \tilde K^T \tilde K$, can be relaxed in the form of the following linear matrix inequality:
\[
\left(
\begin{array}{cc}
W   &  V^T    \\
V & \tilde W  \\ 
\end{array} 
\right) \succeq
\left(
\begin{array}{c}
K^T \\
\tilde K^T  \\ 
\end{array} 
\right)
\left(
\begin{array}{cc}
K & \tilde K  \\ 
\end{array} 
\right)
\quad\Longleftrightarrow\quad
\left( 
\begin{array}{ccc}
W   &  V^T       & K^T        \\
V   & \tilde W & \tilde K^T \\ 
K   & \tilde K &  I \\ 
\end{array} 
\right)\succeq 0,
\]
$W_{ij}, \tilde W_{ij}, V_{ij}\ge 0,\ V_{jj} =0$ for all $i,j$.
\end{remark}

\begin{theorem}\label{theo-3}
Let $G_A$ and $G_B$ be undirected vertex-labeled invertible graphs on $n$ and $m$ vertices without loops, respectively. Assume $G_B$ is arbitrarily bridgeable over the first $k_B$ vertices $\{1,\cdots, k_B\}$. Then 
\[
\Lambda_{HL}(G_C)\le  \Lambda^{opt}_{HL}(G_A,G_B)
\equiv \overline{\Lambda}_{HL}^{sir}(G_A,G_B)
\le \overline{\Lambda}_{HL}^{sdp}(G_A,G_B) 
\le \Lambda_{HL}(G_A),
\]
for any graph $G_C={\mathcal B}_K(G_A,G_B)$ which is constructed from graphs $G_A, G_B$ by bridging the  vertices of $G_A$ to the first $k_B$ vertices of $G_B$ through an $(n,m)$-bipartite graph $G_K$ such that $K_{ij}=0$ for $j=k_B+1, \cdots, m$.
\end{theorem}

\noindent P r o o f. The set 
\[
\{ (K,W),\ 
K_{ij}\in \{0,1\}, \ \  W_{ij}\in\N^+_0, \ \ W_{jj} = \sum_{l} K_{lj}
\quad\hbox{for all}\ i,j, 
\]
\[
\quad \sum_{k,l} K_{kl} \ge 1, W\succeq K^T K \}
\]
of feasible integer matrices $K,W$ for (\ref{homolumosir}) is a subset of the set: 
\[
\{ (K,W),\ 
0\le K_{ij} \le 1,\ \  W_{ij}\ge 0,\ \ W_{jj} = \sum_{l} K_{lj}, 
\quad \hbox{for all}\ i,j,
\] 
\[
\quad \sum_{k,l} K_{kl}\ge 1, W\succeq K^T K
\},
\]
of real matrices $K,W$ that are feasible for (\ref{homolumosdp}). From this fact we conclude the inequality $\overline{\Lambda}_{HL}^{sir}(G_A,G_B) \le \overline{\Lambda}_{HL}^{sdp}(G_A,G_B)$. The inequality $\overline{\Lambda}_{HL}^{sdp}(G_A,G_B) \le \Lambda_{HL}(G_A)$ follows from the fact that 
\[
\left( 
\begin{array}{cc}
I - \mu A^{-1}  &  K  B^{-1} \\
B^{-1} K^T  &  I - \mu B^{-1} + B^{-1}  W  B^{-1}
\end{array}
\right) \succeq 0 \quad \Longrightarrow I - \mu A^{-1}  \succeq 0,
\]
that is $1/\mu \ge \lambda_{max}(A^{-1})$ 
and so $\mu \le \check{\lambda}^+(G_A)$. Similarly, we obtain $I + \eta A^{-1}  \succeq 0$ and, consequently, $\eta \le - \hat{\lambda}^-(G_A)$. Therefore,
$\mu+\eta\le \Lambda_{HL}(G_A)$, as claimed.
\hfill$\diamondsuit$

\medskip

\begin{figure}[htp]
\begin{center}
\includegraphics[width=0.3\textwidth]{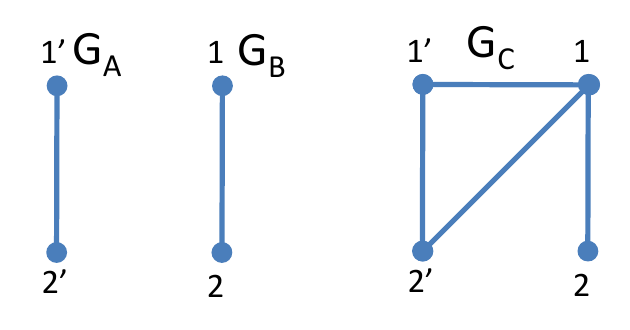}
\end{center}
\caption{
Simple graphs $G_A$ and $G_B$ (left) and the bridged graph $G_C$ with the maximal HOMO-LUMO spectral gap which can be constructed by bridging $G_A$ and $G_B$ over the vertex 1 of $G_B$ ($k_B=1$) to the vertices of $G_A$  (right).
}
\label{fig-quarticfamily}
\end{figure}

\begin{example}
In  Figure~\ref{fig-quarticfamily} (left) we show two simple graphs $G_A$ and  $G_B$ having the spectrum $\sigma(G_A) = \sigma(G_B) = \{1,-1\}$, i.e. $\Lambda_{HL}(G_A)=2$. The graph $G_B$ is arbitrarily bridgeable over the vertex $1$. The optimal bipartite graph $G_K$ bridging $G_B$ to $G_A$ with $k_B=1$ has the adjacency matrix $K=(1,1)^T$. The optimal bridged graph $G_C$ is shown in Figure~\ref{fig-quarticfamily} (right) and it has the spectrum $\sigma(G_C)=\{2.1701, 0.3111, -1, -1.4812\}$, i.e. $\Lambda^{opt}_{HL}(G_A,G_B) = \overline{\Lambda}_{HL}^{sir}(G_A,G_B)= 1.3111$. On the other hand, it turns out that $\overline{\Lambda}_{HL}^{sdp}(G_A,G_B)=1.67597$. Hence we have the strict inequalities 
\[
\Lambda^{opt}_{HL}(G_A,G_B) \equiv  \overline{\Lambda}_{HL}^{sir}(G_A,G_B)
< \overline{\Lambda}_{HL}^{sdp}(G_A,G_B) 
< \Lambda_{HL}(G_A),
\]
in this example. 

\end{example}

\begin{figure}[htp]
\begin{center}
\includegraphics[width=3.5truecm]{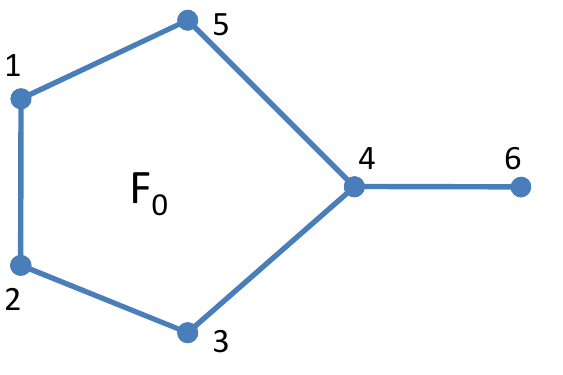}
\ \ 
\includegraphics[width=2truecm]{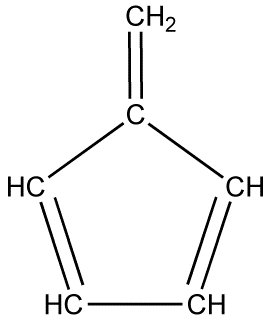}
\end{center}
\caption{
An example of an invertible graph $F_0$ (left) representing the chemical organic molecule of fulvene (right).}
\label{fig-fulvene}
\end{figure}

In Figure~\ref{fig-fulvene} (left) we show the graph $F_0$ on 6 vertices representing the fulvene organic molecule (5-methylidenecyclopenta-1,3-diene) (right). The spectrum consists of the following eigenvalues:
\[
 \sigma(F_0)= \{ 2.1149, 1, 1/q, -0.2541, -q, -1.8608  \},
\]
where $q=(\sqrt{5}+1)/2$ is the golden ratio. The HOMO-LUMO spectral gap $\Lambda_{HL}(F_0)= 0.872134$. It is easy to verify that the graph $G_B\equiv F_0$ is arbitrarily bridgeable over the following subsets of vertices: 
$\{5\}, \{4\}, \{3\}, \{2\}, \{1\}$ for $k_B=1$, $ \{4,5\}, \{2,5\}, \{3,4\}, \{2,4\}, \{1,4\}, \{1,3\}, \{1,2\}$ for $k_B=2$, and $\{2,4,5\}, \{1,3,4\}, \{1,2,4\}$ for $k_B=3$ (cf. Pavl\'{\i}kov\'a and \v{S}ev\v{c}ovi\v{c} \cite{Pavlikova2016}).

\section{Lower bounds for the optimal HOMO-LUMO spectral gap}

In this section, our aim is to derive lower bounds for the optimal HOMO-LUMO separation gap $\Lambda^{opt}_{HL}(G_A,G_B)$. Similarly as in derivation of upper bounds we will construct the lower bound by means of a solution to a certain nonlinear optimization problem. 

The idea is based on construction of upper bounds for the maximal eigenvalues $\lambda_{max}(\pm C^{-1})$ of the inverse matrices $C^{-1}$ and $-C^{-1}$. Here $C$ is the adjacency matrix of the bridged graph $G_C={\mathcal B}_K(G_A, G_B)$. This way we obtain a lower bound for the first positive and negative eigenvalues of $C$ yielding the HOMO-LUMO spectral gap for $G_C$.

The maximal eigenvalue $\lambda_{max}(C^{-1})$ can be expressed by means of the Rayleigh quotient, and, consequently, it can be estimated as follows:
\begin{eqnarray*}
\lambda_{max}(C^{-1}) &=& \max_{\Vert z\Vert^2=1} z^T C^{-1} z 
= (Q z)^T 
\left( 
\begin{array}{cc}
A^{-1} & 0 \\
0      & B^{-1}
\end{array}
\right)
Q z \\
&=& (x-K B^{-1} y)^T A^{-1} (x-K B^{-1} y) + y^T B^{-1} y
\\
&\le& \lambda_{max}(A^{-1}) \Vert  x-K B^{-1} \Vert^2
+ \lambda_{max}(B^{-1}) \Vert  y \Vert^2,
\end{eqnarray*}
where $z=(x,y)\in \R^n\times\R^m$ and the matrix $Q$ is given as in (\ref{invC2}). Analogously,
\[
 \lambda_{max}(-C^{-1})\le \lambda_{max}(-A^{-1}) \Vert  x-K B^{-1} \Vert^2
+ \lambda_{max}(-B^{-1}) \Vert  y \Vert^2.
\]
To estimate the right hand side of the estimate for $\lambda_{max}(\pm C^{-1})$ we apply the following auxiliary lemma proved in \cite{Pavlikova2016}.

\begin{lemma}\cite[Lemma 1]{Pavlikova2016}
\label{lemmaMax}
Assume that $D$ is an $n\times m$ matrix and $\alpha,\beta>0$ are positive constants. Then, for the optimal value $\gamma^*$ of the following constrained optimization problem:
\begin{equation}\label{minD}
\begin{array}{rl}
\gamma^*= \max & \alpha \Vert x- D y\Vert^2 +\beta \Vert y\Vert^2 \\
{\rm s. t.} & \Vert x\Vert^2 + \Vert y\Vert^2  = 1,\ \  x\in\R^n, y\in\R^m,
\end{array} 
\end{equation}
we have the explicit expression:
\begin{eqnarray*}
 \gamma^* 
 &=& \max\left\{ \gamma,\ \frac{(\gamma-\alpha)(\gamma-\beta)}{\alpha\gamma} \in\sigma(D^T D)\right\} \\
 &=&  \frac{\alpha(\omega^* + 1)+\beta + \sqrt{(\alpha(\omega^* +1) +\beta)^2-4\alpha\beta}}{2},
\end{eqnarray*}
where $\omega^*=\max\{\sigma(D^T D)\}$ is the maximal eigenvalue of the matrix $D^T D$. 
\end{lemma}

With help of the previous lemma we obtain the upper estimate:
\[
  \lambda_{max}(\pm C^{-1})\le \frac{\alpha^\pm(\omega^* + 1)+\beta^\pm + \sqrt{(\alpha^\pm(\omega^* +1) +\beta^\pm)^2-4\alpha^\pm\beta^\pm}}{2}
\]
where $\alpha^\pm=\lambda_{max}(\pm A^{-1})$,
$\beta^\pm=\lambda_{max}(\pm B^{-1})$, and, 
\[
\omega^* = \max \sigma(B^{-1}K^T K B^{-1}).
\]
Indeed, for the matrix $D= K B^{-1}$ we have $D^T D = B^{-1} K^T K B^{-1}$. 
The maximal eigenvalue of the matrix $B^{-1}K^T K  B^{-1}$ can be expressed by means of a solution to the semidefinite programming problem:
\begin{eqnarray}
\label{omega}
\omega^* &=& \max \sigma(B^{-1} K^T K  B^{-1}) 
= \min_{B^{-1} K^T K  B^{-1}\preceq \omega I} \omega \nonumber
\\
&=& \min_{\omega} \ \  \omega 
\\
&& s.t.\ \ 
\left( 
\begin{array}{cc}
\omega I     &  B^{-1}  K^T \\
K  B^{-1} &  I
\end{array}
\right) \succeq 0. \nonumber
\end{eqnarray}
Since 
\begin{eqnarray*}
\Lambda_{HL}(G_C)&=&\Lambda_{HL}({\mathcal B}_K(G_A, G_B)) 
\\ 
&=& \check{\lambda}^+(G_C) - \hat{\lambda}^-(G_C) = \frac{1}{\lambda_{max}(C^{-1})} + \frac{1}{\lambda_{max}(-C^{-1})},
\end{eqnarray*}
and the optimal value $\gamma^*$ is an increasing function of $\omega^*$ we obtain the following lower bound $\underline{\Lambda}_{HL}^{sir}(G_A, G_B)\le \Lambda_{HL}^{opt}(G_A, G_B)$ for the optimal HOMO-LUMO spectral gap $\Lambda_{HL}^{opt}(G_A, G_B)$, where

\begin{eqnarray}
\label{lowerSIR}
\underline{\Lambda}_{HL}^{sir}(G_A, G_B) 
&=& 
\frac{2}{\alpha^+(\omega^* + 1)+\beta^+ + \sqrt{(\alpha^+(\omega^* +1) +\beta^+)^2-4\alpha^+\beta^+}} \nonumber \\
&& + \frac{2}{\alpha^-(\omega^* + 1)+\beta^- + \sqrt{(\alpha^-(\omega^* +1) +\beta^-)^2-4\alpha^-\beta^-}}, \nonumber 
\\
&&\nonumber 
\\
\hbox{where} &&  \omega^* =  \min_{\omega, K} \ \  \omega \nonumber
\\
&& s.t.\ \ 
\left( 
\begin{array}{cc}
\omega I     &  B^{-1} K^T \\
K B^{-1} &  I
\end{array}
\right) \succeq 0. \label{lowerconsSIR}
\\
&& K_{i,j}\in\{0,1\}, \ \ \hbox{for each}\ i,j,\ \ \sum_{k,l} K_{kl}\ge 1. 
\nonumber
\end{eqnarray}
Similarly, as in the construction of the upper bound, we can relax the condition $K_{i,j}\in\{0,1\}$ by the box constraint 
\begin{equation}
0\le K_{i,j} \le 1, \ \ \hbox{for each}\ i,j, \label{lowerconsSDP}
\end{equation}
in order to construct the full semidefinite relaxation for the lower bound $\underline{\Lambda}_{HL}^{sdp}(G_A, G_B)$.

\begin{theorem}\label{theo-4}
Let $G_A$ and $G_B$ be undirected vertex-labeled invertible graphs on $n$ and $m$ vertices without loops, respectively. Assume $G_B$ is arbitrarily bridgeable over the first $k_B$ vertices $\{1,\cdots, k_B\}$. Then 
\[
\underline{\Lambda}_{HL}^{sdp}(G_A, G_B)
\le 
\underline{\Lambda}_{HL}^{sir}(G_A, G_B)
\le
\Lambda^{opt}_{HL}(G_A,G_B).
\]
\end{theorem}

\section{Additional constraints imposed on the bridging bipartite graph}

In practical applications one may impose additional constraints on the bridging bipartite graph $G_K$. For example, in computational chemistry the so-called chemical molecules play important role. The structural graph $G$ of a chemical molecule has all vertices of the degree less or equal to 3. If the goal is to construct a bridged graph $G_C={\mathcal B}_K(G_A, G_B)$ representing a chemical molecule with the maximal degree $M_{d}$, we can add additional constraint:

\begin{equation}
\sum_k C_{ik} \le M_{d}, \quad \hbox{for all}\ i,\ \ \hbox{where}
\quad 
C  = \left( 
\begin{array}{cc}
A &  K \\
K^T & B
\end{array}
\right).
\label{constraintChemGraph} 
\end{equation}
The inequality (\ref{constraintChemGraph}) is linear in the $K$ variable and it can be easily added to all nonlinear optimization problems (\ref{homolumoopt}), (\ref{homolumosir}), (\ref{homolumosdp}), (\ref{lowerconsSIR}), (\ref{lowerconsSDP}). The computational results of construction of graphs with the maximal degree $M_d=3$ are presented in the next section. 

Another useful constraint imposed on the bridging graph $G_K$ is the min-max box constraints:
\begin{eqnarray}
\label{constraintMinMax} 
&&\underline{L}^{A}_i\le \sum_k K_{ik} \le \overline{L}^{A}_i, \quad\hbox{for all}\ i=1,\cdots, n,
\\
&&\underline{L}^{B}_j\le \sum_k K_{kj} \le \overline{L}^{B}_j, \quad\hbox{for all}\ j=1,\cdots, k_B,
\end{eqnarray}
representing the box constraints for minimal and maximal number of edges in the bridging graph $G_K$ pointing from the graph $G_A$ to $G_B$. Again, such a box constraint can be easily added to (\ref{homolumoopt}), (\ref{homolumosir}), (\ref{homolumosdp}), (\ref{lowerconsSIR}), (\ref{lowerconsSDP}).

\begin{center}
\end{center}

\section{Computational results}

\begin{figure}[htp]
\begin{center}\includegraphics[width=0.5\textwidth]{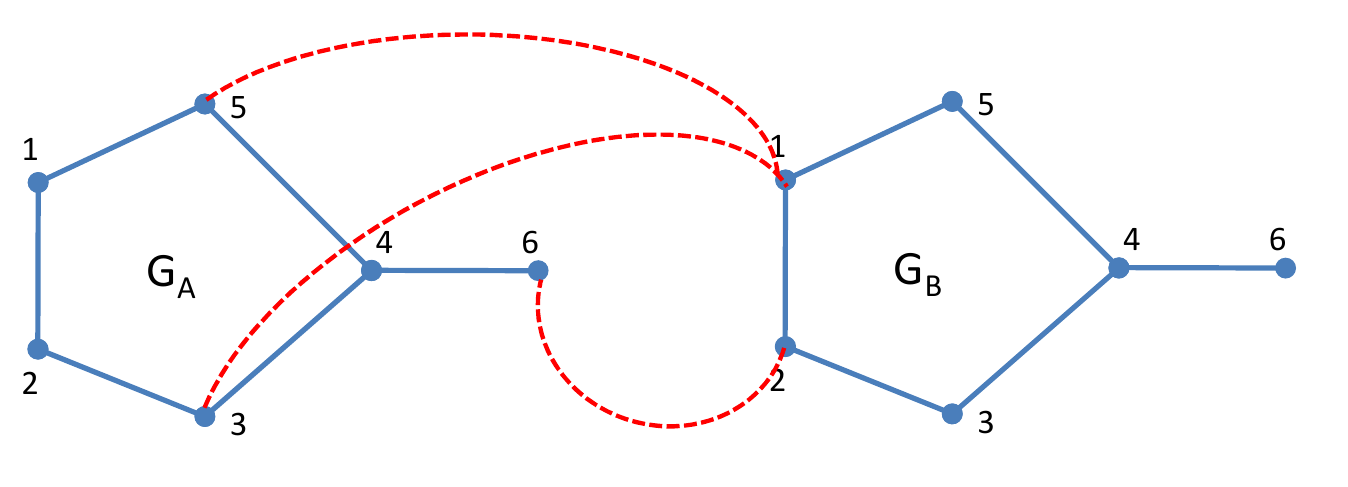}
\\
(a)

\vskip 4truemm
\includegraphics[width=0.5\textwidth]{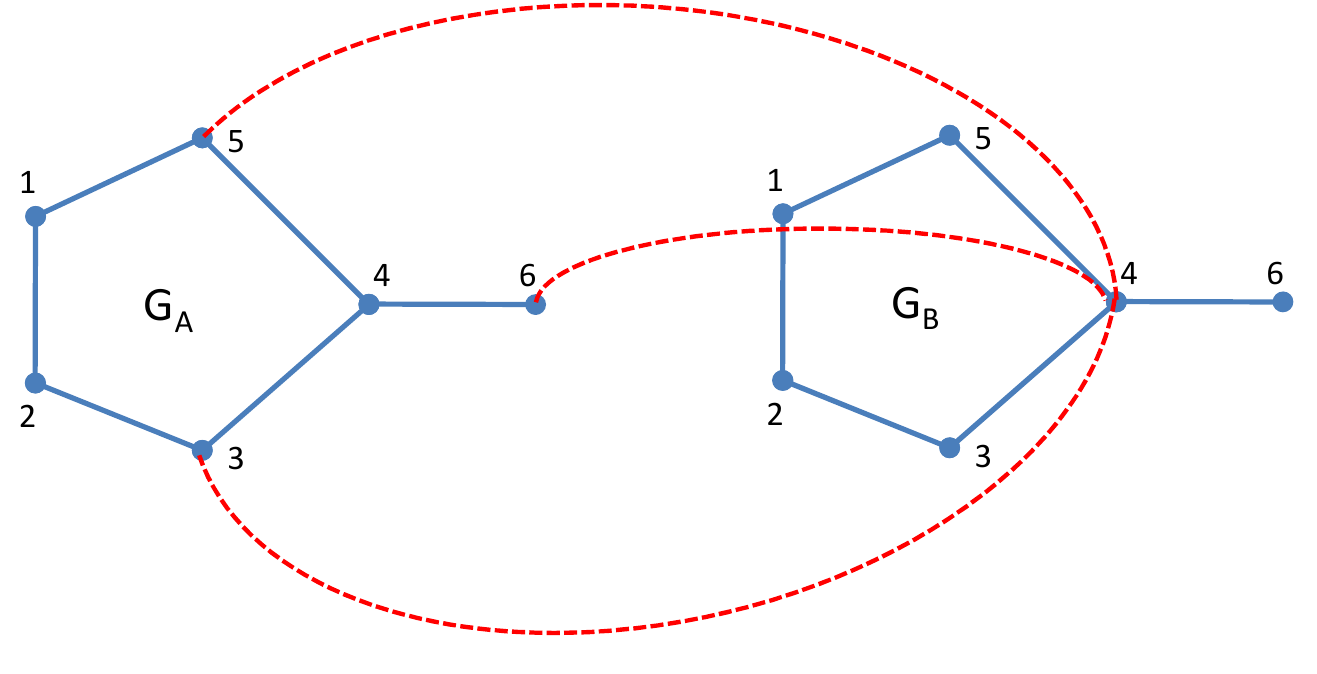}
\\
(b)

\vskip 4truemm
\includegraphics[width=0.75\textwidth]{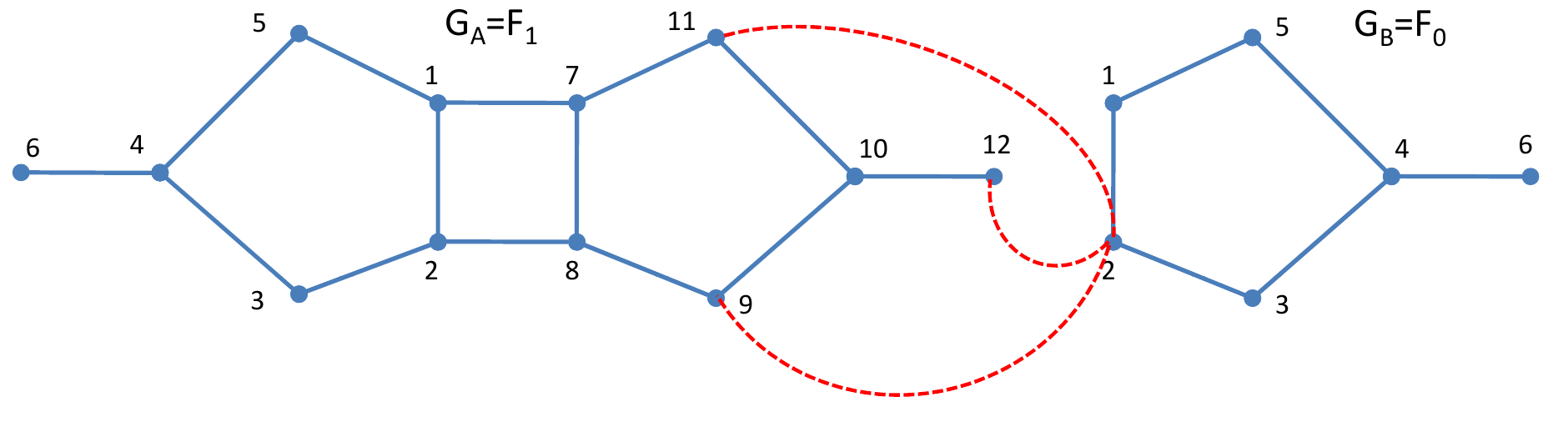}
\\
(c)
\end{center}

\caption{Results of optimal bridging of the fulvene graph $G_B=F_0$ through the vertices $\{1,2\}$ to $G_A=F_0$ a); through the vertices $\{1,4\}$ to $G_A=F_0$ b); and through the vertices  $\{1,2\}$ to $G_A=F_1$ c).}
\label{fig:opt-bridgeF0F0}
\end{figure}

\begin{table}[ht]
\caption{The computational results and comparison of various semidefinite relaxations. The first two columns describe the graph $G_A$ and $G_B$ with the chosen set of bridging vertices. The optimal value $\Lambda_{HL}^{opt}=\overline{\Lambda}_{HL}^{sir}$ is shown in bold in the middle column. The upper $\overline{\Lambda}_{HL}^{sdp}$ and lower bounds $\underline{\Lambda}_{HL}^{sdp}$,  $\underline{\Lambda}_{HL}^{sir}$ are also presented together with computational times in seconds computed on Quad core Intel 1.5GHz CPU with 4 GB of memory.} 
\label{tab-results}

\begin{center}
\scriptsize
\begin{tabular}{c|c||c|c|c|c|c}
$G_A$ & $G_B$ & $\underline{\Lambda}_{HL}^{sdp}$ & $\underline{\Lambda}_{HL}^{sir}$ & $\Lambda_{HL}^{opt}=\overline{\Lambda}_{HL}^{sir}$ & $\overline{\Lambda}_{HL}^{sdp}$ & \hbox{bridging $G_B \mapsto G_A$ } \\
\hline\hline
$F_0$ &  $F_0$ & $0.233688$ & $0.531664$ & $\bf 0.74947$       & $0.87214$ & $1\mapsto 3,5;\ \  2\mapsto 6$ \\
      &$(1,2)$ & $(0.27s)$ & $(3.38s)$ & $(83s)$   & $(2.2s)$  &                            \\
\hline
$F_0$ &  $F_0$ & $0.333126$ & $0.72678$ & $\bf 0.85828$       & $0.87214$ & $1\mapsto \emptyset;\ \  4\mapsto 3,5,6$ \\
      &$(1,4)$ & $(0.31s)$ & $(4.75s)$ & $(36s)$   & $(2.2s)$  &                            \\
\hline
$F_0$ &  $F_0$ & $0.333126$ & $0.719668$ & $\bf 0.81389$       & $0.87214$ & $1\mapsto 4;\ \  3\mapsto 4$ \\
      &$(1,3)$ & $(0.31s)$ & $(4.27s)$ & $(75s)$   & $(2.2s)$  &                            \\
\hline
$F_1$ &  $F_0$ & $0.163626$ & $0.450022$ & $\bf 0.56655$       & $0.56666$ & $1\mapsto \emptyset; \ \ 2\mapsto 9,11,12$ \\
      &$(1,2)$ & $(0.28s)$ & $(7.65s)$ & $(12470s)$   & $(2.2s)$  &                            \\
\hline
$P_4$ &  $P_4$ & $0.472136$ & $0.86953$& $\bf 1.06418$      & $1.23607$ & $2\mapsto 2,4;\ \  3\mapsto 1,3$ \\
      &$(2,3)$ & $(0.27s)$& $(2.18s)$ & $(12.6s)$ & $(2.2s)$  &                            \\
\hline
$P_6$ &  $P_4$ & $0.367365$ & $0.811369$ & $\bf 0.87366$      & $0.89008$ & $1\mapsto 4,6;\ \  3\mapsto 4,6$ \\
      &$(1,3)$ & $(0.26s)$ & $(4.6s)$ & $(59s)$  & $(2.1s)$  &                            \\
\hline
$P_6$ &  $P_4$ & $0.367365$ & $0.737641$ & $\bf 0.87321$      & $0.89008$ & $2\mapsto 4,6;\ \  3\mapsto 1,3$ \\
      &$(2,3)$ & $(0.26s)$ & $(3.41s)$ & $(57s)$  & $(2.1s)$  &                            \\
\hline
$P_{10}$&$P_4$ & $0.252282$ & $0.523808$ & $\bf 0.56837$      & $0.56926$ & $2\mapsto 8,10;\ \  3\mapsto \emptyset$ \\
      &$(2,3)$ & $(0.26s)$  & $(6.32s)$ & $(4109s)$ & $(2.6s)$  &                            \\
\hline
$T_4$ &  $P_4$ & $0.38832$ & $0.73094$& $\bf 0.93258$      & $0.95452$ & $2\mapsto 3,8$ \\
      &$(2)$   & $(0.31s)$& $(1.57s)$ & $(12s)$ & $(2.31s)$  &                            \\
\hline

\end{tabular}

\end{center}

\end{table}

In this section we present computational results. In Table~\ref{tab-results} we present results of construction of the optimal bridging by a bipartite graph for various sets of bridged graphs $G_A$ and $G_B$. First, we chose the fulvene graph $F_0$ as the graph $G_B$ and set $k_B=2$. The graph $G_B\equiv F_0$ is arbitrarily bridgeable through the pairs vertices $\{1,2\}, \{1,3\}, \{1,4\}$ (cf. \cite{Pavlikova2016}). We show the results of the optimal value  $\Lambda_{HL}^{opt}=\overline{\Lambda}_{HL}^{sir}$ for target graphs $G_A = F_0$ and $G_A=F_1$ (see Figure~\ref{fig:opt-bridgeF0F0}). We also presented upper and lower bounds obtained by means of the full semidefinite relaxation. Among the tested examples the maximal HOMO-LUMO gap was attained in the case when $G_B=F_0$ was bridged to $G_A=F_0$ through vertices $\{1,4\}$. Solving mixed integer semidefinite program (\ref{homolumosir}) is time consuming (see Table~\ref{tab-results}). On the other hand, we provided upper and lower bounds which had been obtained efficiently by means of the full semidefinite relaxation technique. A graphical presentation of optimal bridging of fulvene graphs can be seen in Figure~\ref{fig:opt-bridgeF0F0}.

\begin{figure}[htp]
\begin{center}
\includegraphics[width=0.4\textwidth]{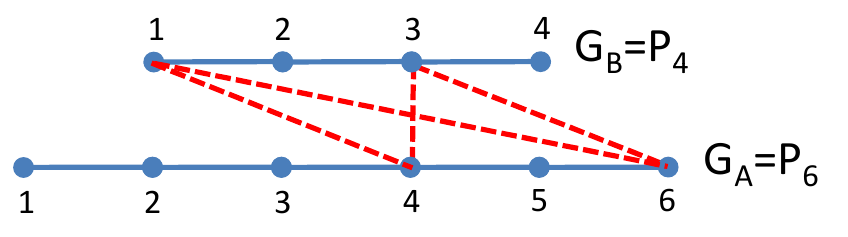}
\qquad 
\includegraphics[width=0.4\textwidth]{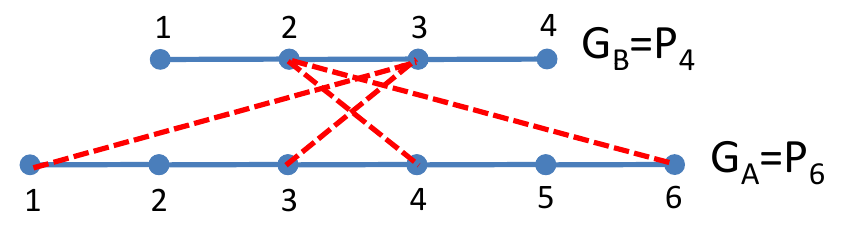}
\\
(a) \hskip 3truecm (b)

\vskip 6truemm
\includegraphics[width=0.5\textwidth]{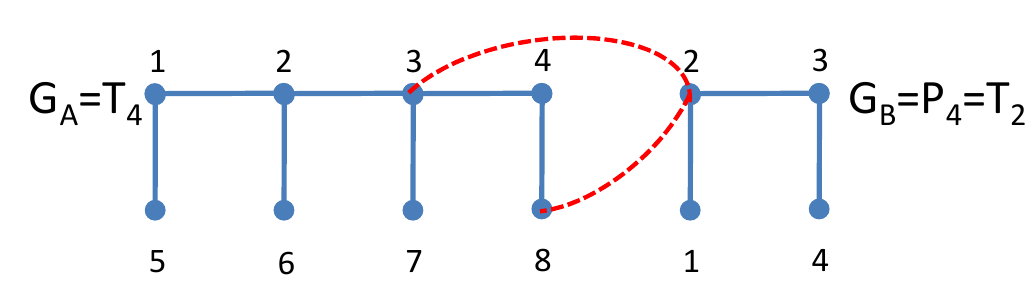}
\\
(c)
\end{center}

\caption{Results of optimal bridging of the graph $G_B=P_4$ through the vertices $\{1,3\}$ to $G_A=P_6$ a); through the vertices $\{2,3\}$ to $G_A=P_6$ b); and through the vertex  $\{2\}$ to $G_A=T_4$ c).}
\label{fig:opt-bridgeP4P6}
\end{figure}

The next set of examples consists of bridging a simple path $G_B=P_m, m=4$ to the path $G_A=P_n, n=4,6$. An illustration of optimal bridging of $P_4$ to $P_6$ over various pairs of vertices is shown in Figure~\ref{fig:opt-bridgeP4P6}. 

The last example is the optimal bridging of $G_B=P_4$ to the graph $G_A=T_{2k}$, where $T_{2k}$ is the graph consisting of the simple path $P_k$ with attached pendant vertices to each vertex of $P_4$. In this case solving the optimal bridging problem yields the bridged graph $G_C$ containing a circle $C_4$ (see  Figure~\ref{fig:opt-bridgeP4P6}, c)).

\begin{figure}[htp]
\begin{center}
\includegraphics[width=0.4\textwidth]{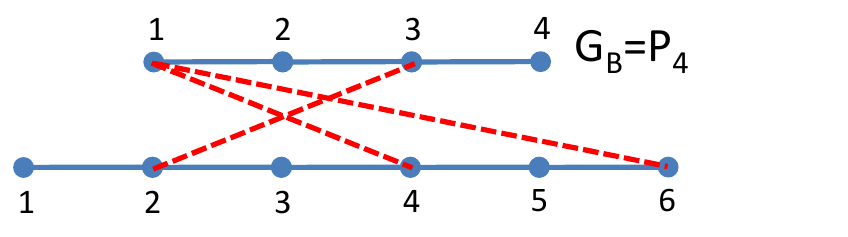}
\qquad 
\includegraphics[width=0.4\textwidth]{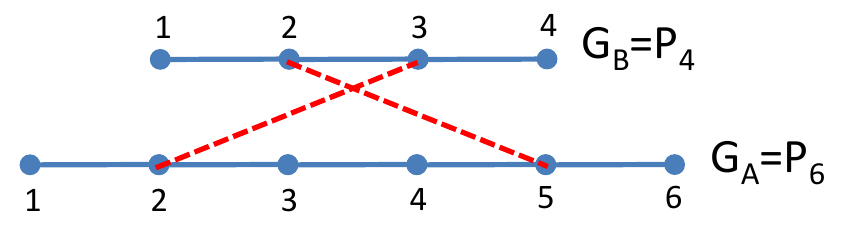}
\\
(a) \hskip 3truecm (b)

\vskip 6truemm
\includegraphics[width=0.5\textwidth]{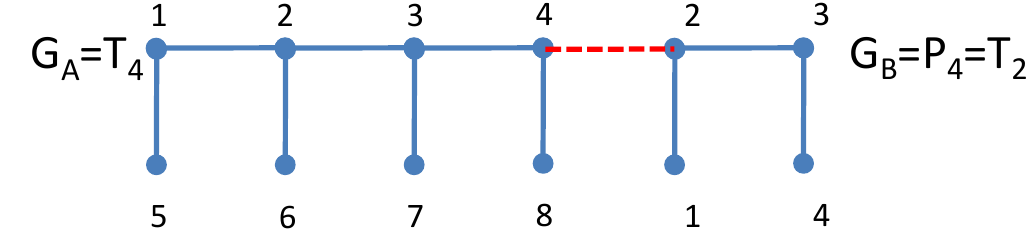}
\\
(c)
\end{center}

\caption{Results of optimal bridging of the graph $G_B=P_4$ through the vertices $\{1,3\}$ to $G_A=P_6$ a); through the vertices $\{2,3\}$ to $G_A=P_6$ b); and through the vertex  $\{2\}$ to $G_A=T_4$ c); with the constraint of maximal degree equal to 3.}
\label{fig:opt-bridgeP4P6-chem}
\end{figure}

In Section 7 we discussed additional constraints imposed on the bridging graph $G_K$. In what follows, we present results of computing the optimal HOMO-LUMO gap and its upper and lower bound under the constraint that the resulting graph $G_C$ represents a chemical molecule with the maximal vertex degree $M_d=3$. The results are summarized in Table~\ref{tab-results-chem} and illustrative examples are shown in Figure~\ref{fig:opt-bridgeP4P6-chem}. In Figure~\ref{fig:opt-bridgeP4P6-chem}, c), we confirmed the well known fact that the comb graph $T_{2k}$ has the maximal HOMO-LUMO gap among all trees on $2k$ vertices with perfect matchings. It was first proved by Kr\v{c} and Pavl\'{\i}kov\'a \cite[Theorem 7]{Pavlikova1990} (see also Zhang and An \cite{Zhang1999}). Interestingly enough, adding additional constraint on maximal degree of vertices considerably reduced computational time for solving the mixed integer semidefinite problem (\ref{homolumosir}). 

\begin{table}

\caption{The computational results and comparison of various relaxations. The chosen graphs and description of columns is the same as in Table~\ref{tab-results}. In this table we present results of optimization when additional constraint of the maximal degree 3 has been imposed.} 
\label{tab-results-chem}
 
\begin{center}
\scriptsize
\begin{tabular}{c|c||c|c|c|c|c}
$G_A$ & $G_B$ & $\underline{\Lambda}_{HL}^{sdp}$ & $\underline{\Lambda}_{HL}^{sir}$ & $\Lambda_{HL}^{opt}=\overline{\Lambda}_{HL}^{sir}$ & $\overline{\Lambda}_{HL}^{sdp}$ & \hbox{bridging $G_B \mapsto G_A$ } \\
\hline\hline
$F_0$ &  $F_0$ & $0.233688$ & $0.507678$ & $\bf 0.720830$       & $0.87214$ & $1\mapsto \emptyset;\ 2\mapsto 6$ \\
      &$(1,2)$ & $(0.31s)$ & $(2.73s)$ & $(7.1s)$   & $(2.9s)$  &                            \\
\hline
$F_0$ &  $F_0$ & $0.233688$ & $0.468053$ & $\bf 0.720830$       & $0.87214$ & $1\mapsto 6; 4\mapsto \emptyset$ \\
      &$(1,4)$ & $(0.31s)$ & $(1.1s)$ & $(2.33s)$   & $(2.85s)$  &                            \\
\hline
$F_0$ &  $F_0$ & $0.333126$ & $0.706635$ & $\bf 0.776875$       & $0.87214$ & $1\mapsto 6; 3\mapsto 6$ \\
      &$(1,3)$ & $(0.35s)$ & $(2.45s)$ & $(8.4s)$   & $(2.82s)$  &                            \\
\hline
$F_1$ &  $F_0$ & $0.163626$ & $0.389941$ & $\bf 0.493727$       & $0.566658$ & $1\mapsto 6; 2\mapsto \emptyset$ \\
      &$(1,2)$ & $(0.38s)$ & $(3.67s)$ & $(13.4s)$   & $(2.83s)$  &                            \\
\hline
$P_4$ &  $P_4$ & $0.472136$ & $0.869530$& $\bf 0.954520$      & $1.23607$ & $3\mapsto \emptyset; 2\mapsto 2$ \\
      &$(2,3)$ & $(0.31s)$& $(1.86s)$ & $(7.8s)$ & $(2.86s)$  &                            \\
\hline
$P_6$ &  $P_4$ & $0.367365$ & $0.811369$ & $\bf 0.828427$      & $0.89008$ & $1\mapsto 4,6; 3\mapsto 2$ \\
      &$(1,3)$ & $(0.36s)$ & $(3.35s)$ & $(22.9s)$  & $(2.83s)$  &                            \\
\hline
$P_6$ &  $P_4$ & $0.367365$ & $0.737641$ & $\bf 0.820751$      & $0.89008$ & $2\mapsto 5; 3\mapsto 2$ \\
      &$(2,3)$ & $(0.33)$ & $(2.73s)$ & $(9.21s)$  & $(2.87s)$  &                            \\
\hline
$P_{10}$&$P_4$ & $0.252282$ & $0.523808$ & $\bf 0.559046$      & $0.56926$ & $2\mapsto \emptyset; 3\mapsto 11$ \\
      &$(2,3)$ & $(0.33s)$  & $(4.78s)$ & $(13.87s)$ & $(2.86s)$  &                            \\
\hline
$T_4$ &  $P_4$ & $0.38832$ & $0.692266$& $\bf 0.890084$      & $0.95452$ & $2\mapsto 4$ \\
      &$(2)$   & $(0.31s)$& $(0.88s)$ & $(1.5s)$ & $(2.11s)$  &                            \\
\hline

\end{tabular}

\end{center}

\end{table}

\section*{Conclusions}

We analyzed spectral properties of graphs which are constructed from two given invertible graphs by bridging them over a bipartite graph. We showed how the HOMO-LUMO spectral gap can be computed by means of a solution to mixed integer semidefinite programming problem. We investigated the optimization problem in which we constructed a bridging graph maximizing the HOMO-LUMO spectral gap. We also provided upper and lower bounds to the optimal value, again expressed as solution to relaxed semidefinite programming problems. Various computational examples were presented in this paper.

\medskip
Received xxxx 20xx; revised xxxx 20xx.
\medskip

\end{document}